\input amstex
\documentstyle{amsppt}
\document
\topmatter
\title
K\"ahler manifolds with homothetic foliation by curves.
\endtitle
\author
W\l odzimierz Jelonek
\endauthor

\abstract{The aim of this paper is to classify compact, simply
connected K\"ahler manifolds which admit totally geodesic,
holomorphic complex homothetic foliations by curves.}
\thanks{MS Classification: 53C55,53C25. Key words and phrases:
K\"ahler manifold, holomorphic foliation, homothetic foliation,
special K\"ahler-Ricci potential, special K\"ahler
potential}\endthanks
 \endabstract
\endtopmatter
\define\G{\Gamma}
\define\DE{\Cal D^{\perp}}
\define\e{\epsilon}
\define\n{\nabla}
\define\om{\omega}
\define\r{\rightarrow}
\define\w{\wedge}
\define\k{\diamondsuit}
\define\th{\theta}
\define\p{\partial}
\define\a{\alpha}

\define\lb{\lambda}

\define\1{D_{\lb}}
\define\2{D_{\mu}}
\define\0{\Omega}

\define\bp{\overline{\partial}}

\define\De{\Cal D}

\define\dl{\delta}

\define\m{(M,g,J)}
\define \E{\Cal E}
\bigskip
{\bf 0. Introduction. } The aim of the present paper is to
classify compact, simply connected K\"ahler manifolds $(M,g,J)$,
$dim M=2n>2$,  admitting a global, complex homothetic foliation
$\Cal F$ by curves which is totally geodesic and holomorphic. A
foliation $\Cal F$ on a Riemannian manifold $(M,g)$ is called
conformal if
$$L_Vg=\a(V)g$$ holds on $T\Cal F^{\perp}$ where $\a$ is a one
form vanishing on $T\Cal F^{\perp}$.  A foliation $\Cal F$ is
called homothetic if is conformal and $d\a=0$  (see [V], [Ch-N]).
Complex homothetic foliations by curves on K\"ahler manifolds were
recently  classified locally in [Ch-N]. Let $\De$ be a
distribution determined by a foliation $\Cal F$.
 By $\E$ we shall denote the $(\dim M-2)$-dimensional distribution
 which is the orthogonal complement of $\De$ in $TM$. If $X$ is a local unit section of $\De$
 then  $\{X,JX\}$
 is a local orthonormal basis of $\De$ and the function $\kappa=
 \sqrt{(div_{\E}X)^2+(div_{\E}JX)^2}$ does not depend on the
 choice of $X$. It turns out that $\kappa=(n-1)|\a|$. We classify  compact, simply connected  K\"ahler
 manifolds admitting a global, complex homothetic foliation $\Cal F$
 satisfying the conditions
 $U\ne \emptyset$ where $U=\{x\in M: |\a|\ne 0\}$, i.e. if $\a$ does not vanish identically.
 First we  show that $\m$
admits a global holomorphic Killing vector field with a Killing
potential, which is a special K\"ahler potential. Next we use
slightly generalized results of Derdzinski and Maschler [D-M-1],
[D-M-2]. Our results rely heavily on the papers [D-M-1], [D-M-2].
As a corollary we prove that every compact, simply connected
K\"ahler manifold admitting a holomorphic, totally geodesic
homothetic foliation with $\a_{x_0}\ne 0$ at least at one point
$x_0\in M$ is a holomorphic $\Bbb{CP}^1$-bundle over  Hodge
manifold, $M=\Bbb P(\Cal L)$ where $\Cal L$ is a holomorphic line
bundle with curvature $\Omega=s\Omega_N$, where $s\ne0$ and
$\Omega_N$ is the K\"ahler form of $(N,h,J)$.  The leaves of the
foliation are the fibers $\Bbb{CP}^1$ of the bundle. If $\a=0$
then $M$ is a product of a Riemannian surface and a K\"ahler
manifold. The result was partially proved in [J] and for the
completeness we cite some of the calculations and proofs from [J]
in our present paper.
\bigskip
{\bf 1. Principal field.} Let $\m$ be a $2n$-dimensional K\"ahler
manifold with a $2$-dimensional $J$-invariant distribution $\De$.
Let $\frak X(M)$ denote the algebra of all differentiable vector
fields on $M$ and $\G(\De)$ denote the set of local sections of
the distribution $\De$. If $X\in\frak X(M)$ then by $X^{\flat}$ we
shall denote the 1-form $\phi\in \frak X^*(M)$ dual to $X$ with
respect to $g$, i.e. $\phi(Y)=X^{\flat}(Y)=g(X,Y)$. By $\0$ we
shall denote the K\"ahler form of $\m$ i.e. $\0(X,Y)=g(JX,Y)$. Let
us denote by $\E$ the distribution $\DE$, which is a
$2(n-1)$-dimensional, $J$-invariant distribution. By $h,m$
respectively we shall denote the tensors $h=g\circ (p_{\De}\times
p_{\De}),m=g\circ (p_{\E}\times p_{\E})$, where $p_{\De},p_{\E}$
are the orthogonal projections on $\De,\E$ respectively. It
follows that $g=h+m$. By $\om$ we shall denote the K\"ahler form
of $\De$ i.e. $\om(X,Y)=h(JX,Y)$ and by $\0_m$ the K\"ahler form
of $\E$ i.e. $\0_m(X,Y)=m(JX,Y)$.    For any local section
$X\in\G(\De)$ we define $div_{\E}X=tr_m\n
X^{\flat}=m^{ij}\n_{e_i}X^{\flat}(e_j)$ where
$\{e_1,e_2,...,e_{2(n-1)}\}$ is any basis of $\E$ and $[m^{ij}]$
is a matrix inverse to $[m_{ij}]$, where $m_{ij}=m(e_i,e_j)$. Note
that if $f\in C^{\infty}(M)$ then $div_{\E}(f X)=fdiv_{\E}X$ in
the case  $X\in \G(\De)$. Let $\xi\in\G(\De)$ be a unit local
section of $\De$. Then $\{\xi,J\xi\}$ is an orthonormal basis of
$\De$. Let $\eta(X)=g(\xi,X)$ and $J\eta=-\eta\circ J$ which means
that $J\eta(X)=g(J\xi,X)$.  Let us denote by $\kappa$ the function
$$\kappa=\sqrt{(div_{\E}\xi)^2+(div_{\E}J\xi)^2}.\tag 1.1$$
The function $\kappa$ does not depend on the choice of a section
$\xi$. In fact, if $\xi'=a\xi+b J\xi$, where $a,b\in
C^{\infty}(dom \xi)$ and $a^2+b^2=1$ is another unit section of
$\De$, then $J\xi'=-b\xi+aJ\xi$ and
$$\gather
(div_{\E}\xi')^2+(div_{\E}J\xi')^2=(a div_{\E}\xi+b
div_{\E}J\xi)^2+(-b div_{\E}\xi+a div_{\E}J\xi)^2\tag 1.2
\\=(div_{\E}\xi)^2+(div_{\E}J\xi)^2.\endgather$$ Hence $\kappa$ is
a well defined, continuous function on $M$, which is smooth in the
open set $U=\{x:\kappa(x)\ne 0\}$. We shall now show that on $U$
there is a smooth, global unit  section $\xi\in \G(U, \De)$
defined uniquely up to  a sign such that $div_{\E}J\xi=0$. Namely,
if $\xi'$ is a local unit section of $\G(U,\De)$ then then the
section $\xi=\frac 1{\kappa}
((div_{\E}\xi')\xi'+(div_{\E}J\xi)J\xi')$ satisfies
$div_{\E}J\xi=\frac1{\kappa}((div_{\E}\xi')(div_{\E}J\xi')-(div_{\E}J\xi')(div_{\E}\xi'))=0$
and does not depend on the choice of $\xi'$. On the other hand it
is clear that the only other such smooth section is $-\xi$. The
section $\xi$ constructed above and defined on $U\subset M$ we
shall call the principal section of $\De$  (see also [G-M]). Note
that $div_{\E}\xi=\kappa$.
\medskip
{\bf 2. Complex homothetic foliations.} We start with  (see [V],
[Ch-N]):

{\bf Definition.} A foliation $\Cal F$ on a Riemannian manifold
$(M,g)$ is called conformal if
$$L_Vg=\a(V)g$$ holds on $T\Cal F^{\perp}$ where $\a$ is a one
form vanishing on $T\Cal F^{\perp}$.  A foliation $\Cal F$ is
called homothetic if is conformal and $d\a=0$.

\medskip
In the rest of the paper we assume that $dim\Cal F=2$ and $\Cal F$
is complex, which means that for an associated distribution $\De$
we have $J\De=\De$. Let us write $\a(X)=g(\zeta,X)$.  Then
$div_{\E}J\zeta=0,div_{\E}\zeta=(n-1)|\a|^2$, which means that in
$U$ the field $\xi=\frac1{|\a|}J\zeta$ is the principal field. Let
$\eta=\xi^{\flat}$.  Note that $\kappa=(n-1)|\a|$.  If we assume
that $\Cal F$ is totally geodesic then $d_{\E}|\a|=0$ i.e.
$d|\a|(X)=0$ if $X\in\E$.  In fact  ([Ch-N]) since $d\a=0$
$Xg(\zeta,\zeta)=2g(\n_X\zeta,\zeta)=2g(\n_{\zeta}\zeta,X)=0$.  A
distribution $\De$ is called holomorphic if $L_XJTM\subset \De$
for any $X\in\G(\De)$. Hence if $\De$ is holomorphic then for any
$X\in\G(\De),Y \in \frak X(M),Z\in\G(\E)$ we have $g(L_XJY,Z)=0$.
Let $A=\n X$. Note that $L_XJ(Y)=A\circ JY-J\circ AY=[A,J]Y$.
Consequently
$$g(AJY,Z)=g(JAY,Z)\tag 2.1$$ for $Y \in \frak X(M),Z\in\G(\E)$.
Let us write $\n_{J\xi}J\xi=p^*\xi$ for a certain function $p^*\in
C^{\infty}( U)$.
\medskip
{\bf Proposition 2.1.}  {\it  Let a foliation $\Cal F$ on a
K\"ahler manifold $\m$ be totally geodesic, holomorphic complex
homothetic foliation by curves.  Then in $U$}
$$\gather d\eta=0,\n_\xi\xi=0,d\ln |\a|=-(|\a|+p^*)\eta,
dp^*\w\eta=0,\\
\n_X\eta(Y)=\frac12|\a|m(X,Y)-p^*J\eta(X)J\eta(Y)\endgather$$
\medskip
{\it Proof.}
 The distribution $\Delta=\{X\in
TU:\eta(X)=0\}$ defined in $U$ is integrable. From (2.1)  it
follows that $\n\eta(JX,JY)=\n\eta(X,Y)$ for $X,Y\in\G(\E)$. If we
take $X=\xi, A=\n\xi$ in (2.1) then we obtain
$$g(AX,Y)+g(X,AY)=\a(\xi)g(X,Y)$$
for any $X,Y\in \G(\E)$. On the other hand since $d\a=0$ we get
$g(AX,Y)=g(X,AY)$ on $\E$. Hence
$$g(AX,Y)=\frac12|\a|g(X,Y)\tag 2.2$$
Consequently, since $\De$ is totally geodesic, we obtain
$$\n_X\eta(Y)=\frac{|\alpha|}{2}m(X,Y)+p\eta(X)J\eta(Y)-p^*J\eta(X)J\eta(Y),\tag
2.3$$ where $\n_\xi\xi=pJ\xi$.  It is also clear that
$$d\eta=p \ \eta\w J\eta\tag 2.4$$ and
$$dJ\eta=|\a|\0- (p^*+|\a|)\eta\w J\eta.\tag 2.5$$
Thus the distribution $\E_{|U}$ is the so called $B_0$
-distribution defined in [G-M]. If dim $M=2n\ge 4$ we also have
 $p=g(\n_{\xi}\xi,J\xi)=0$ and consequently $\n_{\xi}\xi=0,
d\eta=0$.  In fact from (2.4) we get

$$dp \w \eta\w J\eta+ |\a|p\eta\w\0=0,$$
and consequently $p=0$ in $U=\{x\in M:|\a|\ne0\}$. From (2.5) we
obtain
$$ (d|\a|+|\a|(|\a|+p^*)\eta)\w\0=d(p^*+|\a|)\w\eta\w J\eta.\tag 2.6$$
Since $d|\a|\in\bigwedge^1\De$ we get
$(d|\a|+|\a|(|\a|+p^*)\eta)\in\bigwedge^1\De$.  Thus it follows
from (2.6) that  $d|\a|+|\a|(|\a|+p^*)\eta=0$. Hence
$$d\ln |\a|=-(|\a|+p^*)\eta, dp^*\w\eta=0.\tag
2.7$$ Thus
$$\n_X\eta(Y)=\frac12|\a|m(X,Y)-p^*J\eta(X)J\eta(Y).\k\tag
2.8$$

\bigskip
{\bf 3.} {\bf Examples and Killing vector fields with special
K\"ahler potential} First we give a definition
 \medskip {\bf Definition.} {\it A nonconstant
function $\tau\in C^{\infty}(M)$, where $\m$ is a K\"ahler
manifold, is called a special K\"ahler potential if the field $X=
J(\n \tau)$ is a Killing vector field and, at every point with
$d\tau\ne 0$ all nonzero tangent vectors orthogonal to the fields
$X,JX$ are  eigenvectors of  $\n d\tau$}.

Let $\tau$ be a special K\"ahler potential on a K\"ahler manifold
$\m$,
$$\Cal V=span\{\n\tau,J\n\tau\}$$ on $U=\{x\in M:d\tau(x)\ne0\}$ and
let $\Cal F$ be a foliation on $U$ given by the integrable
distribution $\Cal V$. Then $\Cal F$ is a totally geodesic,
holomorphic complex conformal foliation. We have $\a =2\frac MQ
d\tau,\zeta=2\frac MQ \n\tau$, where by $M$ we denote the
eigenvalue of the Hessian $H^{\tau}$ corresponding to the
distribution $\Cal H=\Cal V^{\perp}$  and $Q=g(\n\tau,\n\tau)$.
\medskip

 {\bf Proposition 3.1.}  {\it Let  $X=J(\n \tau)$ be a
  holomorphic  Killing field on
   a K\"ahler manifold  $(M,g,J)$. Then $X$ is an eigenfield
   of the tensor $H^\tau$ if and only if
    $$dQ=2\Lambda d\tau\tag 3.1$$ for a certain function $\Lambda$, which
    is then an eigenfunction of $H^\tau$.}
\medskip
{\it Proof.}  Since $X$ is a  holomorphic   Killing field  then
$$dd^c\tau(Y,Z)=2H^\tau(JY,Z),$$ since $dd^c=2i\p\bp$.  Note that
$L_Xd^c\tau=-L_X(d\tau\circ J)=L_X(d\tau)\circ J=d(L_X\tau)\circ
J=0$. Hence $X\lrcorner dd^c\tau=-d(X\lrcorner d^c\tau)=-d(g(\n
\tau,\n \tau)=-d Q$.
 On the other hand $dd^c\tau(Y,Z)=2H^\tau(JY,Z)$, hence
$-d Q=2H^\tau(JX,.)=-2H^\tau(\n \tau,.)$. The field $\n \tau$ is
an eigenfield of $H^\tau$ if and only if $H^\tau(\n
\tau,.)=\Lambda g(\n \tau,.)=\Lambda d\tau$. It follows that $\n
\tau$ is an eigenfield of $H^\tau$ if $d Q=2\Lambda d\tau$ and
then $H^\tau_{|\Cal V}=\Lambda id_{\Cal V}$.$\k$

\medskip
{ \bf Proposition 3.2.} {\it $\Cal F$ is a totally geodesic,
holomorphic complex homothetic foliation by curves if and only if
the field $X= J(\n \tau)$ is an eigenvector of the Ricci tensor
$\rho$ of $\m$.}
\medskip
{\it Proof.}  We have to show that $d\a=0$ where $\a=2\frac MQ
d\tau$.  Since $dQ=2\Lambda d\tau$ it is equivalent to $dM\w
d\tau=0$.  On the other hand $d\Delta\tau=2\rho(\n\tau,.)$ since
$\n\tau$ is holomorphic and $\Delta\tau=-g(g,\n
d\tau)=-(2\Lambda+2(n-1)M)$. Hence $dM\w d\tau=0$ if and only if
$\rho(\n\tau,.)\w d\tau=0$ which means that $\n\tau$ is an
eigenfield of the Ricci tensor $\rho$.$\k$

\medskip
{\it Remark.}  It follows that compact K\"ahler manifolds
admitting special K\"ahler-Ricci potential described in [D-M-1],
which are holomorphic $\Bbb{CP}^1$ bundles over K\"ahler Einstein
manifolds, give examples of totally geodesic, holomorphic complex
homothetic foliation by curves.  The leaves of the foliation are
the fibers $\Bbb{CP}^1$ of the bundle.

\medskip

If a Killing field $X$ has a special K\"ahler potential then the
    distribution  $\Cal V=span\{X,JX\}$ is totally geodesic.  In
    fact
   $\n_XX=\Lambda JX=-\n_{JX}JX$ i $\n_XJX=J\n_XX=-\Lambda X$.
   Now we prove that if $dim M\ge 6$ then a Killing field $X$ with a special K\"ahler potential
  $X$ is in $U$ an eigenvector of the
   Ricci tensor $S$ of  $\m$.  This fact does not hold if $dim
   M=4$ as Derdzinski shows in [D].

  {\bf Theorem 3.3.}  {\it Let  $X=J(\n \tau)$ be
  a holomorphic Killing field  with a special K\"ahler potential
   on a K\"ahler manifold  $(M,g,J)$ and dim $M\ge 6$. Then

  (a)  $SX=\lb X$ in $U$, where $\lb\in C^{\infty}(U)$ and $  d(\Delta \tau)=2\lb d\tau$,

  (b)  $[S,T]=0$ and $ \n_XS=0$ where $T=\n X$

  (c)  $dQ=2\Lambda d\tau$

  If  $J\circ T_{|\Cal V}=-\Lambda id_{|\Cal V}$ and    $J\circ T_{|\Cal H}=-M id_{|\Cal
  H}$, $Q=g(X,X)$ then  $dQ=2\Lambda d\tau$ and  (a) implies (b). If $\psi(Y)=g(X,Y)=d^c\tau$ then $d\psi(Y,Z)=2g(TY,Z)$.
  In particular
   $\Delta \psi= 2\rho(X,.)= 2\lb \psi$, where $\rho$ is the Ricci tensor.  What is more $dd^c\tau=2\Lambda\om_{\Cal V}+2M\om_{\Cal H }$.}
 \bigskip

{\it Proof.} Note first that $$H^\tau(Y,Z)=g(\n_Y\n
\tau,Z)=-g(J\n_YJ\n \tau,Z)=-g(JTY,Z)\tag 3.2$$ and functions
$\Lambda,M$ are smooth in $U$. We have $2\Lambda+2(n-1)M=-\Delta
\tau$ and $\Lambda=\frac1{Q}H^\tau(X,X)\in C^{\infty}(U)$.
 Hence
$$\gather dQ(Y)=g(\n Q,Y)=-2g(TX,Y)=\\-2g(J\n_X\n
\tau,Y)=2H^\tau(X,JY)=2\Lambda g(X,JY)=\tag3.3\\2\Lambda g(J\n
\tau,JY)=2\Lambda d\tau(Y).\endgather$$

 Note that $JTX=-\Lambda X  $ which implies
 $J\n_YTX+JT^2Y=-Y\Lambda X-\Lambda TY$. Consequently $$
 g(JR(Y,X)X,JY)-g(TY,TY)= -Y\Lambda g(X,JY)-\Lambda g(TY,JY)$$
 and  $$R(Y,X,X,Y)-||TY||^2= Y\Lambda Y\tau+ \Lambda g(JTY,Y)$$  where $R(X,Y,Z,W)=g(R(X,Y)Z,W)$ for any $X,Y,Z,W\in TM$.
 In particular if  $Y\in\Cal H$ then
 $$R(Y,X,X,Y)=M^2||Y||^2-  \Lambda
 M||Y||^2=M(M-\Lambda)||Y||^2.\tag 3.4$$

On the other hand  $2TX=-\n Q$ i $R(Y,X)Y+T^2Y=-\frac12\n_Y\n Q$
which implies
$$  R(Y,X,X,Y)=||TY||^2-\frac12H^{ Q}(Y,Y).$$

Hence
$$  R(Y,X,X,Y')=g(TY,TY')-\frac12H^{ Q}(Y,Y').$$

 For $Y\in\Cal H$  we get $
\frac12H^{ Q}(Y,Y)=\Lambda M||Y||^2$. Hence:

  $$ \gather \lambda Q=2\Lambda^2+(2n-2)M^2+\frac12\Delta Q,\tag
  3.5\\
  \frac12H^{ Q}(X,X)=\Lambda^2 ||X||^2\tag 3.6\\
   \frac12H^{ Q}(JX,JX)=\Lambda^2 ||X||^2-R(JX,X,X,JX),\tag
   3.7\endgather$$

  If $Y\in\Cal H$ then  $$  \frac12H^{ Q}(Y,Y)=\Lambda M||Y||^2\tag 3.8$$
  and
$$R(Y,X,X,Y)=M(M-\Lambda)||Y||^2.\tag 3.9$$

\medskip
{\it Corollary.}   $-\frac12\Delta Q
=-(2\Lambda^2+(2n-2)M\Lambda)+\frac1{ Q}R(JX,X,X,JX)= \lambda
Q-(2\Lambda^2+(2n-2)M^2)$ and  $K(X\w JX)=\frac1{
Q^2}R(JX,X,X,JX)= \lambda+\frac{2(n-1)}{ Q}M(\Lambda-M)$, where
$\lb=\frac1{ Q}\rho(X,X)$.

 \medskip
Now we show that the function $\Lambda,M$ satisfy equations
$d\Lambda\w d\tau=0=dM\w d\tau$ if $dim M\ge 6$ which is
equivalent to the fact that $X$ is an eigenfield of the Ricci
tensor $S$.
  Let
$\om_1=\om^{\Cal V},\om_2=\om^{\Cal H}$, where $\om$ is the
K\"ahler form of $(M,g,J)$. If $\overline{\eta}(Y)=g(X,Y),
\eta(Y)=-g(JX,Y)$ then $\om_1=\frac1{ Q}\eta\w\overline{\eta}$. On
the other hand
$$ d\overline{\eta}=2\Lambda\om_1+2M\om_2,  d\eta=0,
d\om_1=-d\om_2,\n\om_1=-\n\om_2.\tag 3.10$$ Hence
$$d\om_1=-\frac1{ Q} \eta\w
d\overline{\eta}=-2\frac1{ Q}(\Lambda\eta\w\om_1+M\eta\w\om_2).$$
We also have
$$d\Lambda\w\om_1-2(\Lambda-M)\frac1{ Q}(\Lambda\eta\w\om_1+M\eta\w\om_2)+dM\w\om_2=0.$$
Consequently
$$d\Lambda\w\om_1+(dM-2\frac1{ Q}(\Lambda-M)M\eta)\w\om_2=0.$$
Thus   $\n\Lambda\in \G(\Cal V)$ and this relation remains true
also for  $dim M=4$ and $dM=-2\frac1{ Q}(M-\Lambda)M\eta$ which
means that $dM=2\frac1{ Q}(\Lambda-M)Md\tau$.$\k$

 We have  $2\Lambda+2(n-1)M=-\Delta \tau$ and consequently $2d\Lambda+2(n-1)dM=-d\Delta \tau=-2\lambda
 d\tau$, where $\Delta \tau=-tr_gH^\tau$.  Note that  $TX=\Lambda JX$.  Hence $div
 (TX)=tr\{Z\r\n_Z(TX)\}=tr\{R(Z,X)X+T^2Z\}=g(SX,X)-||T||^2=\lambda Q-||T||^2$.
On the other hand  $div(\Lambda JX)=tr\{Z\r Z\Lambda JX+\Lambda
 JTZ\}=d\Lambda(JX)-\Lambda(2\Lambda+2(n-1)M)$. Consequently we
 get

 $$
 \lambda Q-2\Lambda^2-2(n-1)M^2=d\Lambda(JX)-2\Lambda^2-2(n-1)\Lambda
 M.$$  Hence $  d\Lambda(JX)=-2(n-1)M(M-\Lambda)+\lambda Q$
 which means that
 $\lb Q=q Q+2(n-1)M(M-\Lambda)$,
 where  $ d\Lambda=-qd\tau$, and then $ d\Lambda(JX)=-d\Lambda(\n
 \tau)=q Q$.  Hence $q Q=-2(n-1)M(M-\Lambda)+\lambda Q$.
 Consequently
 $  Q d\Lambda=- Q qd\tau= (2(n-1)M(M-\Lambda)-\lambda Q)d\tau$.
 Since    $d\Lambda+(n-1)dM=-\lambda d\tau$ we get  $ (n-1) Q dM
 =(-2(n-1)M(M-\Lambda)+\lambda Q-\lambda Q)
 d\tau=-2(n-1)M(M-\Lambda)d\tau$.  On the other hand
 $q-\lb=\frac{2(n-1)}{ Q}\Lambda(\Lambda-M)$ and we get $q=K(X\w JX)$.

 Hence

    $$\gather  Q dM =2M(\Lambda-M)d\tau,\tag 3.11\\
      Q d\Lambda= (2(n-1)M(M-\Lambda)-\lambda Q)d\tau.\tag
      3.12\endgather$$

  There exists a constant $c$   such that
  $\frac{ Q}M=2(\tau-c)$  or $M=0$ in $U$. In fact if $V=\{x\in
  U:M(x)\ne0\}\ne \emptyset$ then in $V$
 $d\frac{ Q}M= \frac{d Q}{M}- Q\frac{dM}{M^2}=2\frac{\Lambda}M d\tau +2(1-\frac{\Lambda}M)d\tau=2d\tau$
  which implies $ d(\frac{ Q}M-2\tau)=0$ and consequently $V=U$.     Since the set $M'=\{x:X(x)\ne
   0\}$ is connected we obtain

   $$\frac{ Q}M=2(\tau-c).\tag 3.13$$

  Hence  if $\a$ is not identically $0$ then
  $$\a=d\ln |\tau-c|.\tag 3.14$$   We can assume that $c=0$
  replacing $\tau$ by $\tau-c$.
  In the case dim$M=4$  we have to assume, that
   $X$ is an eigenfield of the Ricci tensor $S$ to obtain the above relations.

Let  $Y\in \Cal H$, then  $d\tau(Y)=d\tau(JY)=0$. Hence for any
$Z\in\frak X(M)$ we get

$$  \n_Zd\tau(Y)+d\tau(\n_ZY)=0,  \n_Zd\tau(JY)+d\tau(J\n_ZY)=0$$
which implies
$$\gather \n_ZY^{\Cal V}=-\frac1{ Q}(H^\tau(Z,Y)\n \tau-H^\tau(Z,JY)J\n
\tau)=\\
-M\frac1{ Q}(g(Z,Y)\n \tau+\om(Z,Y)J\n \tau).\endgather$$

 If we assume that $Z\in \Cal H$ then  $[Z,Y]^{\Cal V}=2\frac M{ Q}\om(Y,Z)J\n \tau$.

Now it is easy to see just as in [D-M-1] that
$$QR(Y,Z)\n\tau=2(\Lambda-M)M\om(Y,Z)X,\tag 3.15$$
where we assume that $R(u,v)w=\n_u\n_vw-\n_v\n_uw-\n_{[u,v]}w$ (a
different notation then in [D-M-1]).

\medskip

  {\bf Proposition 3.4.}  {\it Let  $X=J(\n \tau)$ be
  a holomorphic Killing field on a  K\"ahler manifold  $(M,g,J)$.
   Then in $U=\{x:X_x\ne0\}$  the following conditions are equivalent:

(a) $X$ is an eigenfield of $H^\tau$

(b)   $v=\n \tau$ is pregeodesic i.e. $\n _vv=\Lambda v$

(c) The distribution $\Cal V=span(X,JX)$ is totally geodesic

(d)The distribution $\Cal V=span(X,JX)$ is an eigendistribution
of}
 $H^\tau$,

 (e)$d Q=2\Lambda d\tau$.
\medskip
{\it Proof.} $(a)\Leftrightarrow (b)$ $X$ is an eigenfield of
$H^\tau$ if and only if $d Q=2\Lambda d\tau$.  We have $\n
Q=-2\n_XX=-2\n_vv$ which means that $\n Q=2\Lambda\n \tau$ if and
only if $\n_vv=\Lambda v$.

$(b)\Rightarrow (c)$ Since $\n_XX=\n_vv=\Lambda v$ then
$\n_Xv=-J\n_XX=-\Lambda Jv$ and $[X,JX]=0$, and we are done.

$(c)\Rightarrow (a)$  Since $Xg(X,X)=0$ then $g(\n_XX,X)=0$ and
$\n_XX=\lambda JX$ which means that $\n_vv=-\lambda v$.

The equivalence $(a)\Leftrightarrow(d)$ is obvious.$\k$

If $d Q=2\Lambda d\tau$ then $\Delta Q=-2q Q+2\Lambda\Delta \tau$.
In fact $\n_X\n Q=2X\Lambda\n \tau+2\Lambda\n_X\n \tau$ and
$$ H^{ Q}(X,Y)=-2qd\tau(X)d\tau(Y)+2\Lambda H^\tau(X,Y).$$

In particular   $H^\tau_{|\Cal V}=Mg$ wtw. gdy  $H^{ Q}_{|\Cal
V}=2M\Lambda g$. Note that $\eta=d^c\tau$ and
$$dd^c\tau=2\frac{\Lambda}Qd\tau\w d^c\tau+2M\om^{\Cal H}.$$
Consequently if we denote $\th=\frac s{2Q}g(X,.)=\frac
s{2Q}d^c\tau$ where $s\in\Bbb R$, then  $d\th=\frac
s{2|\tau-c|}\om^{\Cal H}$.

\bigskip
{\bf 4.   Local holomorphic Killing vector field on $U$.} Let $\m$
be a K\"ahler  manifold of dimension $2n\ge 4$ admitting a global,
complex homothetic foliation $\Cal F$ by curves, which is totally
geodesic and holomorphic . We shall show in this section using the
ideas from [J] that for every $x\in U$ there exists an open
neighborhood $V\subset U$ of $x$ and a function $f\in
C^{\infty}(V)$ such that $X_V=fJ\xi$ is a Killing vector field in
$V$, which we shall call a local special Killing vector field (see
also [Ch-N]). Let $V$ be a geodesically convex neighborhood of $x$
in $U$. Then $V$ is contractible. Note that the form
$\phi=-p^*\eta$, where $p^*=g(\n_{J\xi}J\xi,\xi)$ is closed in
$U$, since by (2.7), $d\phi=-dp^*\w\eta=0$. Consequently there
exists a function $F\in C^{\infty}(V)$ such that
$$dF=\phi=-p^*\eta.\tag 4.1$$
Let $f=\exp\circ F$. From (2.8) it follows that
$$\n_XJ\eta(Y)=\frac{|\a|}{2}\0_m(X,Y)+p^*J\eta(X)\eta(Y),\tag
4.2$$

Now let $\psi=(fJ\xi)^{\flat}=fJ\eta$. We shall show that
$\n_X\psi(Y)=-\n_Y\psi(X)$ which means that the field $fJ\xi$ is a
Killing vector field in $V$.We get

$$\n_X(fJ\eta)(Y)=XfJ\eta(Y)+f\frac{|\a|}{2}\0_m(X,Y)+fp^*J\eta(X)\eta(Y),\tag
4.3$$ Since $Xf=fXF=-fp^*\eta(X)$ we obtain
$$\gather \n_X(fJ\eta)(Y)=f\frac{|\a|}{2}\0_m(X,Y)-fp^*\eta\w J\eta(X,Y)=\tag
4.4\\f\frac{|\a|}{2}\0_m(X,Y)-fp^*\om(X,Y),\endgather$$ which
proves our claim. Note that if $F_1$ is another solution of
$(4.1)$ then $F_1=F+D$ for a certain constant $D\in\Bbb R$.  It
follows that $f_1=\exp F_1=CF$, where $C=\exp D$. Consequently
$X_1=f_1J\xi=CfJ\xi=CX$. Recall here the well known general fact,
that if $X,Y\in\frak{iso}(M)$ are Killing vector fields on
connected  Riemannian manifold $M$, $X\ne 0$ and $Y=fX$ for a
certain $f\in C^{\infty}(M)$ then $f$ is constant.

Let $\phi=f\xi^{\flat}$. Now it is clear that
$$\n_X\phi(Y)=\n_X(f\eta)(Y)=f\frac{|\a|}{2}m(X,Y)-fp^*h(X,Y)\tag 4.5$$
and $\n_X(f\eta)(Y)=\n_Y(f\eta)(X)$. Consequently $d(f\eta)=0$.
It follows that there exists a function $\tau\in C^{\infty}(V)$
such that $f\xi=\n \tau$. Consequently $X_V=fJ\xi=J(\n \tau)$,
which means that $X_V$ is a holomorphic Killing vector field with
a K\"ahler potential $\tau$.
\medskip
{\bf 5. Special Jacobi fields along geodesics in $U$.}  Let
$c:[0,l]\rightarrow M$  be a unit geodesic such that
$c([0,l))\subset U$ and $c(l)\in K=\{x\in M:\kappa(x)=0\}$. A
vector field $C$, which is a Jacobi field along $c$ i.e. $\n_{\dot
c}^2C-R(\dot c,C)\dot c=0$, will be called a special Jacobi field
if there exists an open, geodesically convex neighborhood $V$ of
$c(0)$ such that $C(0)=X_V(c(0)),\n_{\dot c}C(0)=\n_{\dot
c}X_V(c(0))$. If $im\ c\cap V=c([0,\e))$ then it follows that
$X_V(c(t))=C(t)$ for all $t\in[0,\e)$. We have the following
lemma:
\medskip
{\bf Lemma 5.1.} {\it Let us assume that a vector  field $C$ along
a geodesic $c$ is a special Jacobi field along $c$. Then}

$$\lim_{t\rightarrow l}|C(t)|=0,$$
{\it and} $$g(\dot c,C)=0.$$

{\it Proof.} Let us note that $g(\n_{\dot c}C,\dot c)=0$ since
this property is valid for Killing vector fields. It follows that
the function  $g(\dot c,C)$ is constant. Let $k\in (0,l)$. Then
$c([0,k])\subset U$. For every  $t\in [0,k)$ there exists a
geodesically convex open neighborhood $V_t$ of the point $c(t)$
and a special Killing vector field $X_{V_t}=f_tJ\xi$ on $V_t$
defined in Section 4. The field $X_{V_t}$ is defined uniquely up
to a constant factor. From the cover $\{V_t\}:t\in [0,k]$ of the
compact set $c([0,k])$ we can choose a finite subcover $\{
V_{t_1},V_{t_2},...,V_{t_m}\}$. Let $c_i$ be the part of geodesic
$c$ contained in $V_i=V_{t_i}$, i.e. $im\ c\cap V_i=im\
c_i=c((t_i,t_{i+1}))$. We define the Killing vector field $X_i$ in
every $V_i=V_{t_i}$ by induction in such a way that $X_1=X_V$ on
$V_1\cap V$ and $X_i=X_{i+1}$ on $V_i\cap V_{i+1}$. Let
$X_i=f_iJ\xi$. Note that $C(t)=X_i\circ c(t)=f_iJ\xi\circ c(t)$
for $t\in (t_i,t_{i+1})$. Consequently, on $V_i$, $|C|=f_i$.  From
(2.7) and (4.1) it follows that
$$d\ln\kappa=d\ln f_i-\frac{\kappa}{n-1} \eta.$$
Hence
$$\frac d{dt}\ln\kappa\circ c(t)=\frac d{dt}\ln |C(t)|-\frac{\kappa}{n-1}  \eta(\dot c(t)),\tag 5.1$$
and

$$\frac d{dt}\ln\frac{\kappa\circ c(t)}{ |C(t)|}=-\frac{\kappa}{n-1}  \eta(\dot c(t)).\tag 5.2$$
Consequently
$$\ln\frac{\kappa\circ c(k)}{ |C(k)|}-\ln\frac{\kappa\circ c(0)}{ |C(0)|}=-\frac1{n-1}\int_0^k\kappa \eta(\dot
c(t))dt.\tag 5.3$$ Hence
$$\ln|C(k)|=\ln\kappa\circ c(k)-\ln\frac{\kappa\circ c(0)}{ |C(0)|}+\frac1{n-1}\int_0^k\kappa \eta(\dot
c(t))dt.\tag 5.4$$ Note that
$$|\int_0^k\kappa \eta(\dot
c(t))dt|\le \int_0^k|\kappa \eta(\dot c(t))|dt\le \int_0^k\kappa
|\dot c(t)|dt\le \int_0^k\kappa dt.\tag 5.5$$ Let
$\kappa_0=sup\{\kappa(x):x\in c([0,l])\}$. From (5.4) it follows
that
$$\ln|C(k)|\le\ln\kappa\circ c(k)-\ln\frac{\kappa\circ c(0)}{ |C(0)|}+\frac1{n-1}\kappa_0l.\tag 5.6$$
Consequently
$$\limsup_{k\rightarrow l-}\ln|C(k)|\le\lim_{k\rightarrow l-}\ln\kappa\circ
c(k)-\ln\frac{\kappa\circ c(0)}{
|C(0)|}+\frac1{n-1}\kappa_0l=-\infty.\tag 5.7$$ From (5.7) it is
clear that $\lim_{k\rightarrow l-}|C(k)|=0.$ Since $|g(\dot
c(t),C(t))\le |C(t)|$ and $g(\dot c,C)$ is constant it follows
that $|g(\dot c(t),C(t))|\le \lim_{t\rightarrow l}|C(t)|=0$ which
means that $g(\dot c,C)=0$.$\k$
\medskip
{\bf 6. Global holomorphic Killing  vector field on $M$.} From now
on  we  assume that $\m$ is a complete  K\"ahler manifold with
totally geodesic, holomorphic, complex homothetic foliation $\Cal
F$, $\dim M\ge 4$ and the set $U=\{x\in M:\kappa(x)\ne 0\}$ is
non-empty. Let $K=\{x\in M:\kappa(x)=0\}$.
\medskip
{\bf Theorem 6.1.} {\it The set $U$ is connected and the set $K$
has an empty interior.}
\medskip
{\it Proof.} Let $U_1$ be a non-empty component of the set
$U=M-K$. Let $x_0\in U_1$ and let us assume the set $ int\ K\cup
(U-U_1)$ is non-empty. Let $x_1\in int\ K\cup(U-U_1)$. Then
$x_1=\exp_{x_0}lX$ for a certain unit vector $X\in T_{x_0}M$ and
$l>0$. Let $V\subset U_1$ be a geodesically convex neighborhood of
$x_0$ and $X_V=f_VJ\xi$ the local Killing vector field on $V$. Let
$C$ be the Jacobi vector field along the geodesic
$c(t)=\exp_{x_0}tX$ satisfying the initial
conditions:$C(0)=X_V(x_0),\n_{\dot c}C(0)=\n_{\dot c}X_V(x_0)$. It
follows that the field $C$ is a special Jacobi field along $c$. In
particular $g(X,X_V(x_0))=0$ since the geodesic $c$ meets $K$.
Since the set $int\ K\cup(U-U_1)$ is open there exists an open
neighborhood $W\subset int\ K\cup(U-U_1)$ of the point $x_1$. The
mapping $T_{x_0}M\ni Y\rightarrow \exp_{x_0}lY\in M$ is continuous
hence there exists an open neighborhood $P$ of $X$ in $T_{x_0}M$
such that $\exp_{x_0}lP\subset int\ K\cup(U-U_1)$. We can find a
vector $X_1\in P$ such that $g(X_1,X_V(x_0))\ne 0$. The field
$C_1$ along the geodesic $d(t)=\exp_{x_0}tX_1$ defined by the
initial conditions $C_1(0)=X_{V}(x_0),\n_{\dot d}C_1(0)=\n_{\dot
d(0)}X_V$ is a special Jacobi field along a geodesic $d$ which
meets $K$. It follows that $g(\dot d,C_1)=0$ along $d$. In
particular for $t=0$ we obtain $g(X_1,X_V(x_0))=0$ which is a
contradiction. Consequently $int\ K\cup(U-U_1)=\emptyset$.$\k$

\medskip
Let $x_0\in U$ and let $V=\exp_{x_0}W$ be a geodesically convex
open neighborhood of $x_0$, where $W\subset T_{x_0}M$ is a star
shaped open neighborhood of $0\in T_{x_0}M$. For every $X\in W$
let  $l(X)=\sup\{t:tX\in W\}$. Hence if the sphere $S_{\e}=\{X\in
T_{x_0}M:|X|=\e\}$ is contained in $W-\exp_{x_0}^{-1}(K)$ then
$W=\{tX:X\in S_{\e},t\in[0,l(X))\}$. In an open neighborhood
$V'\subset V-K$ of $x_0$ there is defined a special Killing vector
field $X_{V'}=f_{V'}J\xi$. Let $Z=X_{V'}(x_0)\in T_{x_0}M$ so that
$H=\{X\in T_{x_0}M:g_{x_0}(X,Z)=0\}$ is a hyperplane in
$T_{x_0}M$. Note that $\exp_{x_0}^{-1}(V\cap K)\subset H$. Let a
function $k:S_{\e}\rightarrow \Bbb R$ be defined as follows:
$k(X)=l(X)$ if $x\notin H$ and $k(X)=\inf\{t>0:\exp(tX)\in K\}$ if
$X\in H$. The set $W'=\{tX:X\in S_{\e},t\in[0,k(X))\}$ is open and
star shaped. Note that $V''=\exp_{x_0}W'\subset V-K\subset U$,
$V''$ is dense in $V$,  and $V''$ is contractible. It follows that
in $V''$ there is defined a special local  Killing vector field
$X_{V''}$. We can assume that $X_{V'}=X_{V''}$ on $V'$. Now we can
prove:

{\bf Lemma 6.2} {\it On every geodesically convex open set $V$ in
$M$ can be defined a  holomorphic  Killing vector field $X$ such
that for every open geodesically convex set $W\subset V\cap U$ the
restriction $X_{|W}$ is a special Killing vector field on $W$.}
\medskip
{\it Proof.} We shall use the notation introduced above. Since
$V''$ is contractible there exists a special Killing vector field
$X_{V''}$ defined on $V''$. Let us define a differentiable field
$X$ on $V$ by the formula: $X(\exp_{x_0}u)=J_u(1)$ where $u\in W$
and $J_u$ is a Jacobi vector field along a geodesic
$c(t)=\exp_{x_0}(tu)$ satisfying the initial
conditions:$J_u(0)=X_{V''}(x_0),\n_{\dot c}J_u(0)=\n_{u}X_{V''}$.
It is clear that $X$ is a differentiable vector field and that
$X_{|V''}=X_{V''}$. Since the set $V''$ is dense in $V$ it follows
that $X$ is a Killing vector field in $V$. In fact if we write
$T=\n X$ then for every $Y,Z\in\frak X(V)$ we have
$g(TY,Z)=-g(Y,TZ)$ on $V''$ and both sides are differentiable
functions on $V$. Thus the relation remains valid on $V$ which
means that $X$ is a Killing vector field on $V$. Note that $X=0$
on $V\cap K$.  Since  equation (4.5) is valid on the open, dense
subset $V-K$ of $V$ it follows that the form $\phi=-(JX)^{\flat}$
satisfies the relation $\n_Y\phi(Z)=\n_Z\phi(Y)$ for every
$Y,Z\in\frak X(V)$. Consequently $JX= -\n \tau$ for a certain
function $\tau\in C^{\infty}(V)$ and $X=J\n \tau$ on $V$. Since
$\n_{JY}\phi(JZ)=\n_Y\phi(Z)$ on $V''$ and hence on $V$ for every
$Y,Z\in \frak X(V)$ it is clear that $X$ is holomorphic. $\k$
\medskip
Let $V_1,V_2$ be two open, geodesically convex sets. Then the set
$V_1\cap V_2-K$ is connected. The proof of this is similar to the
proof of Theorem 6.1. Since the sets $V_1''\subset V_1-K$,
$V_2''\subset V_2-K$ are contractible there exist special Killing
vector fields $X_1=f_1J\xi, X_2=f_2J\xi$ on these sets. These
Killing fields can be extended on $V_1-K,V_2-K$ in such a way that
$f_i=\exp F_i$ where $F_i$ satisfy   equation (4.1). Consequently
$d(F_1-F_2)=0$ on the connected set $V_1\cap V_2-K$. Hence
$F_1=F_2+D$ in $V_1\cap V_2-K$ for a constant $D\in \Bbb R$.
Consequently $X_1=C_{12}X_2$ where $C_{12}=\exp D$. The fields
$X_1,X_2$ can be extended to the Killing fields on $V_1,V_2$
respectively such that $X_{i|K\cap V_i}=0$. It is clear that for
these extensions which we also denote by $X_1,X_2$ the equation
$$X_1=C_{12}X_2$$
holds on $V_1\cap V_2$.

\medskip
Now we shall prove:
\bigskip
{\bf Theorem 6.3.  }{\it   Let $\m$ be a complete   K\"ahler
manifold of dimension $2n\ge 4$ with totally geodesic, holomorphic
complex homothetic foliations by curves. Let $\a_{x_0}\ne0$ at
least at one point $x_0\in M$. If $H^1(M,\Bbb R)=0$ then there
exists on $M$ a non-zero holomorphic Killing vector field $X=J(\n
\tau)$ with a special K\"ahler potential $\tau$ such that $X$ is
an eigenfield of the Ricci tensor of $\m$.}
\medskip
{\it Proof.} Let $\{V_i\}:i\in I$ be a cover of $M$ by
geodesically convex, open sets $V_i$. Let $X_i$ be a Killing
vector field on $V_i$ constructed in Lemma 6.2. Let $V_i\cap
V_j\ne\emptyset$. Then there exist constants $C_{ij}>0$ such that
$X_i=C_{ij}X_j$ on $V_i\cap V_j$. These constants satisfy the
co-cycle condition $C_{ij}C_{jk}C_{ki}=1$. Consequently the
constants $D_{ij}=\ln C_{ij}\in \Bbb R$ satisfy the co-cycle
condition $D_{ij}+D_{jk}+D_{ki}=0$. It follows that $\{D_{ij}\}$
is a co-cycle in the first \v{C}ech cohomology group
$\breve{H}^1(\{V_i\},\Bbb R)$. Since $\{V_i\}$ is a good cover of
$M$ it follows that $\breve{H}^1(\{V_i\},\Bbb
R)=\breve{H}^1(M,\Bbb R)=0$. Consequently there exists a co-chain
$\{D_i\}\in Z^0(\{V_i\},\Bbb R)$ such that
$\{D_{ij}\}=\dl(\{D_i\})$. This means that $D_{ij}=D_j-D_i$. Let
$C_i=\exp D_i$. Then $C_{ij}=\frac{C_j}{C_i}$. Let us define the
field $X$ on $M$ by the formula $X_{|V_i}=C_iX_i$. Then it is
clear that $X$ is a well defined, global vector field and
$X\in\frak X(M)$. Since $X_{|V_i}$ is a Killing vector field on
every $V_i$ it follows that $X$ is a Killing vector field. Now let
$\phi=-(JX)^{\flat}$. Then $d\phi=0$, since this equation is
satisfied on every $V_i$. On the other hand the first de Rham
group of $M$ vanishes: $H^1(M,\Bbb R)=\breve{H}^1(M,\Bbb R)=0$.
Consequently there exists a function $\tau\in C^{\infty}(M)$, such
that $\phi=d\tau$. Note also that $\n d\phi$ is Hermitian, which
means that $X$ is  holomorphic.   Thus $X=J(\n \tau)$ is a
holomorphic Killing vector field with a Killing potential $\tau$.
Note that in view of  (4.5) the special Killing field constructed
by us is a Killing vector field with a special K\"ahler potential
$\tau$, which is an eigenfield of the Ricci tensor by Prop. 3.1.
$\k$
\medskip
\bigskip
{\bf Corollary  6.4.  }{\it   Let $\m$ be a complete  K\"ahler
manifold of dimension $2n\ge 4$ with totally geodesic, holomorphic
complex homothetic foliations by curves. Let $\a_{x_0}\ne0$ at
least at one point $x_0\in M$ and let $(\tilde M,\tilde g)$ be the
Riemannian universal covering space of $\m$. Then there exists on
$(\tilde M,\tilde g)$ a non-zero holomorphic Killing vector field
$X$ with a special K\"ahler potential such that $X$ is an
eigenfield of the Ricci tensor of $(\tilde M,\tilde g)$.}

\medskip
{\it Remark.}  Note that if $dim M\ge 6$ then for every special
K\"ahler potential $\tau$ the vector field $J(\n\tau)$ is
automatically  eigenfield of the Ricci tensor $S$ of $(M,g)$ (see
Th.3.3). If $dim M=4$ then from Th.6.3 it follows that in our case
$\tau$ is a special K\"ahler-Ricci potential ([D-M-1]).
\bigskip
{\bf 7. Construction of  K\"ahler manifolds.} In our construction
we shall follow B\'erard Bergery (see [Ber], [S]) rather then
Derdzi\'nski and Maschler, although we shall use the
generalization of classification theorem by Derdzi\'nski and
Maschler ([D-M-1]) to classify  K\"ahler manifolds with complex
homothetic foliation by curves. These two approaches are
equivalent. Let $(N,h,J)$ be a simply connected Hodge
manifold,i.e. $(N,h,J)$ is a K\"ahler manifold and the cohomology
class $\{\frac s{2\pi}\0_N\}$ is an integral class. We also assume
that  $\dim N=2m>2$. Let $s\ge0,L>0,s\in \Bbb Q,L\in \Bbb R$ and
$r:[0,L]\rightarrow \Bbb R$ be a positive, smooth function on
$[0,L]$ with $r'(t)>0$ for $t\in (0,L)$, which is even at $0$ and
$L$, i.e. there exists an $\e>0$ and  even, smooth functions
$r_1,r_2:(-\e,\e)\rightarrow \Bbb R$ such that $r(t)=r_1(t)$ for
$t\in[0,\e)$ and $r(t)=r_2(L-t)$ for $t\in(L-\e,L]$. If $s\ne0$
then it is clear that the function $f=\frac2srr'$ is positive on
$(0,L)$ and $f(0)=f(L)=0$. Let $P$ be a circle bundle over $N$
classified by the integral cohomology class $\{\frac
s{2\pi}\0_N\}$ where $\0_N$ is the K\"ahler form of $(N,h,J)$. On
the bundle $p:P\rightarrow N$ there exists a connection form $\th$
such that $d\th=sp^*\0_N$ where $p:P\rightarrow N$ is the bundle
projection. Let us consider the manifold $(0,L)\times P$ with the
metric
$$g=dt^2+f(t)^2\th^2+r(t)^2p^*h,\tag 7.1$$ if $s\ne0$.
 The metric $7.1$ is
K\"ahler if and only if $f=\frac{2rr'}s$. We shall prove it in
section 9. The proof coincides with the proof in [J]. This time we
do not assume that $(N,h,J)$ is Einstein. It is known that the
metric (7.1) extends to a metric on a sphere bundle
$M=P\times_{S^1}\Bbb{CP}^1$ if and only if the function $r$ is
positive and smooth on $(0,L)$, even at the points $0,L$, the
function $f$ is positive, smooth and odd at the points $0,L$ and
additionally $$f'(0)=1,f'(L)=-1.\tag 7.3$$ If $f=\frac{2rr'}s$ for
$r$ as above, then (7.3) means that
$$2r(0)r''(0)=s,2r(L)r''(L)=-s.\tag 7.4$$

{\bf 8.  Circle bundles.  } Let $(N,h,J)$ be a Hodge manifold with
integral class  $\{\frac s{2\pi}\0_N\}$, where $s\in\Bbb Q$ and
let $p:P\rightarrow N$ be a circle bundle with a connection form
$\th$ such that $d\th=s\0_N$ (see [K]). Let us assume that
dim$N=2m=2(n-1)$. Let us consider a Riemannian metric $g$ on $P$
given by
$$g=a^2\th\otimes\th+b^2p^*h\tag 8.1$$
where $a,b\in \Bbb R$. Let $\xi$ be the fundamental vector field
of the action of $S^1$ on $P$ i.e. $\th(\xi)=1, L_{\xi}g=0$. It
follows that $\xi\in\frak{iso}(P)$ and $a^2\th=g(\xi,.)$.
Consequently
$$a^2d\th(X,Y)=2g(TX,Y)\tag 8.2$$
for every $X,Y\in\frak X(P)$ where $TX=\n_X\xi$. Note that
$g(\xi,\xi)=a^2$ is constant, hence $T\xi=0$. On the other hand
$d\th(X,Y)=sp^*\0_N(X,Y)=sh(Jp(X),p(Y))$.  Note that there exists
a tensor field $\tilde J$ on $P$ such that $\tilde J\xi=0$ and
$\tilde J(X)=(JX_*)^H$ where $X=X_*^H\in TP$ is the horizontal
lift of $X_*\in TN$ (i.e. $\th(X_*^H)=0$) and $X_*=p(X)$. Indeed
$L_{\xi}T=0$ and $T\xi=0$ hence $T$ is the horizontal lift of the
tensor $\tilde T$. Now $\tilde J=\frac{2b^2}{sa^2}\tilde T$. Since
$T\xi=0$ we get $\n T(X,\xi)+T^2X=0$ and $R(X,\xi)\xi=-T^2X$. Thus
$g(R(X,\xi)\xi,X)=||TX||^2$ and
$$\rho(\xi,\xi)=||T||^2=\frac{sa^4}{4b^4}2m.$$ Consequently
$$\lb=\rho(\frac{\xi}a,\frac{\xi}a)=\frac1{a^2}||T||^2=\frac{s^2a^2}{2b^4}m.\tag 8.3$$

We shall compute the O'Neill tensor $A$ (see [ON]) of the
Riemannian submersion $p:(P,g)\rightarrow (N,b^2h)$. We have
$$A_EF=\Cal V(\n_{\Cal H E}\Cal H F)+\Cal H(\n_{\Cal H E}\Cal V
F).$$

Let us write $u=\Cal V(\n_{\Cal H E}\Cal H F)$ and $v=\Cal
H(\n_{\Cal H E}\Cal V F)$.  The vertical component of a field $E$
equals $\th(E)\xi$. If $X,Y\in\Cal H$ then
$$g(\n_XY,\frac1a\xi)=\frac1a
(Xg(Y,\xi)-g(Y,\n_X\xi))=-\frac1ag(TX,Y)=\frac1ag(X,TY).\tag 8.4$$
Hence
$u=\frac1ag(E-\th(E)\xi,T(F-\th(F)\xi)\frac{\xi}a=\frac1{a^2}g(E,TF)\xi.$
Note that $\Cal H(\n_Xf\xi)=f\Cal H(\n_X\xi)=fTX$, hence $$v=\Cal
H(\n_{\Cal H E}\Cal V F)=\th(F)T(E)=\frac1{a^2}g(\xi,F)TE.$$
Consequently
$$A_EF=\frac1{a^2}(g(E,TF)\xi+g(\xi,F)TE).\tag 7.5$$
If $U,V\in \Cal H$ then
$$||A_UV||^2=\frac1{a^2}g(E,TF)^2=\frac{s^2a^2}{4b^4}g(E,\tilde J
F)^2.$$ If $E$ is horizontal and $F$ is vertical then
$$A_EF=\frac1{a^2}g(\xi,F)TE.\tag 8.6$$
Hence $A_E\xi=TE$ and $||A_E\xi||^2=||TE||^2=\frac{s^2a^4}{4b^4}$.
It follows that
$$K(P_{E\xi})=\frac{s^2a^2}{4b^4},$$
where $K(P_{EF})$ denotes the sectional curvature of the plane
generated by vectors $E,F$.  If $E,F\in\Cal H$ then
$$K(P_{EF})=K_*(P_{E_*F_*})-\frac{3g(E,TF)^2}{a^2||E\w F||^2},$$
where $E_*$ denotes the projection of $E$ on $M$ i.e. $E_*=p(E)$.
Thus
$$K(P_{EF})=\frac1{b^2}K_0(P_{E_*F_*})-\frac{3s^2a^2g(E,\tilde JF)^2}{4b^4||E\w F||^2},\tag 8.7$$
where $K_0$ stands for the sectional curvature of the metric $h$
on $N$. Applying this we get for any $E\in\Cal H$ the formula for
the Ricci tensor $\rho$ of $(M,g)$:
$$\rho(E,E)=\frac1{b^2}\rho_0(bE_*,bE_*)-\frac{s^2a^2}{2b^4},\tag 8.8$$
where $\rho_0$ is a Ricci tensor of $(M,h)$.  Now we shall find a
formula for $R(X,\xi)Y$ where $X,Y\in \Cal H$. We have
$R(X,\xi)Y=\n T(X,Y)$ and $$\gather \n
T(X,Y)=\n_X(T(Y))-T(\n_XY)=\n^*_{X_*}(\tilde TY^*)+\frac12\Cal
V[X,TY]\tag 8.9\\-(\tilde T(\n^*_{X_*}Y_*))^*=\frac12\Cal
V[X,TY]=-\frac12sp^*\0_N(X,TY)\xi=-\frac{s^2a^2}{4b^2}h(X_*,Y_*)\xi\endgather$$
Consequently $R(X,Y,Z,\xi)=0$ for $X,Y,Z\in\Cal H$,and
$$R(X,\xi,Y,\xi)=-\frac{s^2a^4}{4b^2}h(X_*,Y_*).\tag 8.10$$
Hence $\rho(\xi,X)=0$ if $g(X,\xi)=0$ and $\xi$ is an eigenfield
of $\rho$. Note that the Ricci tensor $\rho$ has at a point $x\in
P $ $k+1$ eigenvalues where $k$ is the number of eigenvalues of
$\rho_0$ at a point $p(x)$.

{\bf 9.   Riemannian submersion $  p:(0,L)\times P\rightarrow
(0,L)$.}  In this case the O'Neill tensor $A=0$. We shall compute
the O'Neill tensor $T$ (see [ON]).  Note that, since
$L_{\xi}\th=0$,  the field $\xi$ is the Killing vector field for
$((0,L)\times P,g)$, where $g=dt^2+f(t)^2\th^2+r(t)^2p^*h$. We
show that $(0,L)\times P$ is a K\"ahler manifold and $\xi$ is a
Killing vector field with the special K\"ahler potential
$\frac{r^2}s$. Let us denote by $Y^*$ the horizontal lift of the
vector $Y\in TN$ with respect to the Riemannian submersion
$p_N:P\rightarrow N$ i.e. $p_N(Y^*)=Y,g(Y^*,\xi)=0$. Let $H=\frac
d{dt}$ be the horizontal vector field for this submersion and
$\De$ be the distribution spanned by the vector fields $H,\xi$. If
$U,V\in \Cal V$ (here $\Cal V$ temporary  denotes the vertical
distribution of the above Riemannian submersion) and $g(U,V)=0$
then $T(U,V)=0$. Let $U\in \Cal V,g(U,\xi)=0$ and $U=U_*^*$ with
$h(U_*,U_*)=1$, then the following formula holds
$$T(U,U)=-rr'H.\tag 9.1$$
In fact $2g(\n_UV,H)=-Hg(U,V)=-2rr'h(U_*,V_*)$ if $U=V$ or $0$ if
$g(U,V)=0$.
 We also have
$$T(\xi,\xi)=-ff'H.\tag 9.2$$
Now we shall prove that the almost complex structure defined by
$$JH=\frac1f\xi,JX=(J_*X_*)^*\text{   for  } X=(X_*)^*\in\E=\De^{\perp}$$
where $X_*\in TN$, is a K\"ahler structure with respect to the
metric $g$. The proof is similar to that in [J] although now we do
not assume that $(N,h)$ is Einstein.  We have for horizontal lifts
$X,Y\in\frak X(P)\subset\frak X((0, L)\times P)$ of the fields
$X_*,Y_*\in\frak X(N)$ ( with respect to the submersion described
in the Section 7):
$$\gather \n
J(Y,X)=\n_Y(JX)-J(\n_YX)=\n^*_{Y_*}(J_*(X_*))^*-\frac12d\th(Y,JX)\xi\\
+T(Y,JX)-J(\n^*_{Y_*}(X_*)^*-\frac12d\th(Y,X)\xi+T(Y,X))=\\
-\frac12sh(JY,JX)fJH-rr'h(Y,JX)H-\frac12sh(JY,X)fH+h(X,Y)rr'JH=0
\endgather$$
if and only if $f=\frac{2rr'}s$. Note that for $X\in \Cal
H=\De^{\perp}$ we have $\n_X\xi=\frac{sf^2}{2r^2}JX$. It follows
that $\De$ is a totally geodesic foliation. In fact $\n_HH=0,
[\xi,H]=0$ and $g(\n_X\xi,H)=g(\n_X\xi,\xi)=0$ for $X\in\Cal H$.
Since the distribution $\De$ is totally geodesic and
two-dimensional  it is clear that $\n J(X,Y)=0$ if
$X,Y\in\G(\De)$.  Now we shall show that
$$\n J(H,X)=\n J(X,H)=0\text {  for   } X\in \E.$$
It is easy to show that $\n_XH=\n_HX=\frac{r'}rX$ and
$$\n_X(JH)=\n_X(\frac1f\xi)=\frac1fT(X)=\frac{sf}{2r^2}JX.$$
On the other hand
$$\n_X(JH)=\n J(X,H)+J(\n_XH)=\n J(X,H)+\frac {r'}rJ(X).$$
Thus $\n J(X,H)=0$ if $f=\frac{2rr'}s$.  Similarly
$$\n_H(JX)=\n_{JX}H=\frac{r'}rJX=\n J(H,X)+J(\n_HX)=\n
J(H,X)+\frac{r'}rJX$$ and $\n J(H,X)=0$. Note that the K\"ahler
form $\0=fdt\w\th+r^2p^*\0_N$ of almost Hermitian manifold
$((0,L)\times P,g,J)$ is closed, which means that the structure
$J$ is almost K\"ahler. Thus $\n J(JX,Y)=-J\n J(X,Y)$ and
consequently $\n_{\xi}J=0$ which finishes the proof. Let
$U,V,W\in\Cal V$ and $g(U,\xi)=g(V,\xi)=g(W,\xi)=0$. Then for
$R(X,Y,Z,W)=g(R(X,Y)Z,W)$ we get
$$R(U,V,\xi,W)=\hat
R(U,V,\xi,W)+g(T(U,\xi),T(V,W))-g(T(V,\xi),T(U,W))=0.$$ On the
other hand since $\De$ is totally geodesic we obtain
$$ R(H,V,\xi,H)=0$$ and consequently $\rho(V,\xi)=0$ which means
that $\De$ is an eigendistribution of the Ricci tensor $\rho$ of
$((0,L)\times P,g,J)$.

From O'Neill formulae it follows also that
$$R(JH,U,V,JH)=0$$ if $g(U,V)=0$ and
$$R(JH,U,U,JH)=\frac{s^2f^2}{4r^4}-\frac{f'r'}{fr},\tag 9.3$$ for
a unit vector field $U$ as above. Note also that the distribution
$\De$ spanned by the vector fields $\xi,H$ is totally geodesic.
Consequently if $X,Y,Z\in\G(\De)$ and $V$ is as above then
$$R(X,Y,Z,V)=0.\tag 9.4$$
On the other hand for $U,V$ as above and with
$g(U,V)=0,g(U,U)=g(V,V)=1$ we get
$$R(U,V,V,U)=\hat
R(U,V,V,U)-g(T(V,V),T(U,U))=\hat R(U,V,V,U)-\frac{(r')^2}{r^2}$$
and consequently
$$\rho(V,V)=\hat\rho(V,V)-(2n-3)\frac{(r')^2}{r^2}g(V,V)+(\frac{s^2f^2}{4r^4}-\frac{f'r'}{fr})g(V,V)$$
for any $V\in\E$.  It means that at any point  the number of
eigenvalues of $\rho$ is $k+1$ where $k$ is the number of
eigenvalues of the Ricci tensor of $(N,h)$ at the corresponding
point. Hence in general the special K\"ahler potential is not the
special K\"ahler-Ricci potential.  However the distributions
$\De,\E$ are stil orthogonal with respect to $\rho$.
\medskip
{\bf Theorem 9.1.}  {\it Let $\Cal F$ be a holomorphic, complex,
homothetic foliation by curves on a simply connected, compact
K\"ahler manifold $\m$. Let us assume that the form $\a$ does not
vanish identically on $M$ and the leaves of $\Cal F$ are totally
geodesic. Then $M=\Bbb P(\Cal L)$ where $p:\Cal L\r N$ is a
holomorphic line bundle over compact, simply connected  Hodge
manifold $(N,h)$, whose curvature equals $\0=s\om^h, s\ne 0$ and
with a metric defined on the dense open subset $M'=(0,L)\times
P\subset M$
$$g=dt^2+(\frac{2rr'}s)^2\th^2+r(t)^2p^*h,$$
where $r$  satisfies the boundary condition described in section
7.  The leaves of $\Cal F$ are the fibers $\Bbb{CP}^1$ of the
bundle $p:\Bbb P(\Cal L)\r N$.  If $\a=0$ then $M=\Bbb{CP}^1\times
N$ where $N$ is simply connected K\"ahler manifold with the
product metric and the leaves of $\Cal F$ are
$\Bbb{CP}^1\times\{y_0\}$ where $y_0\in N$. }

\medskip
{\it Proof.} In view of Th. 6.3 the distribution $\De$ associated
with the foliation $\Cal F$  coincides in an open, dense subset
$U$ of $M$ with $\De=span\{\n\tau,J\n\tau\}$ where $\tau$ is a
special K\"ahler potential. If dim$M=4$ the potential $\tau$ is a
special K\"ahler-Ricci potential. Using the results in section 3
we can apply the methods and proofs from [D-M-1], [D-M-2]. These
proofs are also valid if we assume only that $\tau$ is a special
K\"ahler potential instead a special K\"ahler-Ricci potential and
apply the results from section 3. Hence we do not assume that
$(N,h)$ is an Einstein manifold if $dim M\ge6$. Note that the
function $\tau$ has two critical submanifolds $N,N^*$ of complex
co-dimension $1$, since otherwise $\Cal F$ defined in $M'=M-(N\cup
N^*) $ would not extend to the foliation on the whole of $M$.  In
our case $s\ne 0$ where the curvature form of the bundle $\Cal L$
is $\0=s\om^h$. Hence $\frac s{2\pi}\om^h$ is an integral form and
$(N,h)$ is a Hodge manifold. In the notation of [D-M-1], [D-M-2]
our notation can be translated as follows. Let us denote $u=\frac
s2\xi$ and $Q=\frac{s^2}4f^2=r^2(r')^2$. We also have $\Cal V=\De,
\Cal H=\Cal E$. Then $u$ is a Killing vector field
$$u=J(\n \frac{r^2}2)$$ and $\tau=\frac{r^2}2$, $c=0, a=\frac s2$.  We have
$JH=\frac1f\xi$ and $\xi=\frac1sJ(\n r^2)$. It is easy to see that
$\n_H\xi=f'JH$ and if $X\in\Cal V$ then $\n_X\xi= \frac2s\Lambda
JX=f'JX$, where $\Lambda= \frac s2 f'$ is an eigenvalue of
$H^{\tau}$ corresponding to the eigendistribution $\Cal V$. We
also have $\th=\frac s{2Q}g(u,)=\frac s{2Q}d^c\tau$ and $M=(r')^2$
is the eigenvalue corresponding to the eigendistribution $\Cal H$.
 The distance between $N$ and $N^*$ is
$$L=\int_{\tau_{min}}^{\tau_{max}}\frac{d\tau}{\sqrt{Q}}$$
 Note that $div_{\E}H=2(n-1)\frac{r'}r$ and
$div_{\E}\xi=0$. In particular $\kappa=2(n-1)\frac{r'}r\ne 0$ on
an open and dense subset and $|\a|=2\frac{r'}r$, $\a=2d\ln r$. On
the other hand for every  Hodge manifold $(N,h)$ we can construct
on the manifold $\Bbb P(\Cal L)$ many K\"ahler metrics $g$ in such
a way that fibers of the bundle $p:\Bbb P(\Cal L)\r N$ form a
totally geodesic, holomorphic complex homothetic foliation. If
$\Cal F$ is a holomorphic, complex, homothetic foliation by curves
on a simply connected, complete K\"ahler manifold $\m$  with
$\a=0$ then $M$ is a product of Riemannian surface $\Sigma$ and a
K\"ahler manifold $N$.  This  follows easily from de Rham theorem
and the fact that in this case both distributions $\De,\E$ are
totally geodesic. If $M$ is compact, simply connected then clearly
$\Sigma=\Bbb{CP}^1$. $\k$

\bigskip
\centerline{\bf References.}

\par
\medskip
\cite{Ber}  L. B\'erard Bergery,{\it Sur de nouvelles vari\'et\'es
riemanniennes d'Einstein}, Pu\-bl. de l'Institute E. Cartan
(Nancy) {\bf 4},(1982), 1-60.
\par
\medskip
\cite{Ch-N} S.G.Chiossi and P-A. Nagy {\it Complex homothetic
foliations on K\"ahler manifolds}  Bull. London Math. Soc. 44
(2012) 113-124.
\par
\medskip
\cite{D} A. Derdzi\'nski {\it Killing potentials with Geodesic
Gradients on K\"ahler Surfaces}  arxiv:1104.4132v1 [math.DG] 20
Apr 2011
\par
\medskip
\cite{D-M-1} A. Derdzi\'nski, G. Maschler {\it Special
K\"ahler-Ricci potentials on compact K\" ahler manifolds}, J.
reine angew. Math. 593 (2006), 73-116.
 \par
\medskip
\cite{D-M-2}  A. Derdzi\'nski, G.  Maschler {\it Local
classification of conformally-Einstein  K\"ahler metrics in higher
dimensions}, Proc. London Math. Soc. (3) 87 (2003), no. 3,
779-819.

\par
\medskip
\cite{G-M} G.Ganchev, V. Mihova {\it K\"ahler manifolds of
quasi-constant holomorphic sectional curvatures}, Cent. Eur. J.
Math. 6(1),(2008), 43-75.

\par
\medskip
\cite{J} W. Jelonek {\it K\"ahler manifolds with quasi-constant
holomorphic curvature} Ann. Glob. Anal. and Geom, vol.36,
(2009),143-159.
\par
\medskip
\cite{K} S. Kobayashi {\it Principal fibre bundles with the
1-dimensional toroidal group} T\^ohoku Math.J. {\bf 8},(1956)
29-45.

\par
\medskip
\cite{ON} B. O'Neill, {\it The fundamental equations of a
submersion}, Mich. Math. J. {\bf 13},(1966), 459-469.
\par
\medskip
\cite{S} P. Sentenac ,{\it Construction d'une m\'etrique
d'Einstein  sur la somme de deux projectifs complexes de dimension
2}, G\'eom\'etrie riemannienne en dimension 4 ( S\'eminaire Arthur
Besse 1978-1979) Cedic-Fernand Nathan, Paris (1981), pp. 292-307.

\par
\medskip
\cite{V} I. Vaisman {\it Conformal foliations} Kodai J. Math.
2(1979), 26-37.
\par
\medskip
Institute of Mathematics

Cracow University of Technology

Warszawska 24

31-155      Krak\'ow,  POLAND.

 E-mail address: wjelon\@pk.edu.pl
\bigskip

\enddocument